\magnification=1200

%%%%%%%%%%%%%%%%%%%%%%%%%%%%%%%%%%%%%%%%%%%%%%%%%%%%%%%%%%%%%%%%%%%%%
% Some fundamental parameters

\hsize=125mm             % horizontal size of text
\vsize=195mm             % vertical size of text
\parskip=0pt plus 1pt    % some variable space between paragraphs
\clubpenalty=10000       % no single lines at the end of a page
\widowpenalty=10000      % no single lines at the beginnig of a page
\frenchspacing           % no double spaces after puntuation
\parindent=8mm           % dimension of indentation; affects \item,
                         % \itemitem and normal paragraphs
%%%%%%%%%%%%%%%%%%%%%%%%%%%%%%%%%%%%%%%%%%%%%%%%%%%%%%%%%%%%%%%%%%%%%
% Fonts & Families
%%%%%%%%%%%%%%%%%%%%%%%%%%%%%%%%%%%%%%%%%%%%%%%%%%%%%%%%%%%%%%%%%%%%%
\let\txtf=\textfont
\let\scrf=\scriptfont
\let\sscf=\scriptscriptfont

% \fam0 Roman % Defined in Plain-TeX
\font\frtnrm =cmr12 at 14pt
\font\tenrm  =cmr10
\font\ninerm =cmr9
\font\sevenrm=cmr7
\font\fiverm =cmr5

\txtf0=\tenrm
\scrf0=\sevenrm
\sscf0=\fiverm

\def\rm{\fam0 \tenrm}

% \fam1 Math-Italic % Defined in Plain-TeX
\font\frtnmi =cmmi12 at 14pt
\font\tenmi  =cmmi10
\font\ninemi =cmmi9
\font\sevenmi=cmmi7
\font\fivemi =cmmi5

\txtf1=\tenmi
\scrf1=\sevenmi
\sscf1=\fivemi

 \def\oldstyle{\fam1 \tenmi}

% \fam2 Symbols % Defined in Plain-TeX
\font\tensy  =cmsy10
\font\ninesy =cmsy9
\font\sevensy=cmsy7
\font\fivesy =cmsy5

\txtf2=\tensy
\scrf2=\sevensy
\sscf2=\fivesy

\def\cal{\fam2 }

% \fam4 Times-Italic % Defined in Plain-TeX
\font\tenit  =cmti10
\font\nineit =cmti9
\font\sevenit=cmti7
\font\fiveit =cmti7 at 5pt

\txtf\itfam=\tenit
\scrf\itfam=\sevenit
\sscf\itfam=\fiveit

\def\it{\fam\itfam\tenit}

%\fam 6 Bold % Defined in Plain-TeX
\font\tenbf  =cmb10
\font\ninebf =cmb10 at  9pt
\font\sevenbf=cmb10 at  7pt
\font\fivebf =cmb10 at  5pt

\txtf\bffam=\tenbf
\scrf\bffam=\sevenbf
\sscf\bffam=\fivebf

\def\bf{\fam\bffam\tenbf}

% \fam B AMS-Font 2 % New family
\newfam\msbfam       
\font\tenmsb  =msbm10
\font\sevenmsb=msbm7
\font\fivemsb =msbm5

\txtf\msbfam=\tenmsb
\scrf\msbfam=\sevenmsb
\sscf\msbfam=\fivemsb

\def\msb{\fam\msbfam\tenmsb}
\def\Bbb#1{{\msb #1}}

\def\CC{{\Bbb C}}
\def\ZZ{{\Bbb Z}}
\def\RR{{\Bbb R}}
\def\NN{{\Bbb N}}
\def\QQ{{\Bbb Q}}
\def\PP{{\Bbb P}}

\def\FP{{\cal F}(P)}

\def\th{\tilde{h}}

% \fam 8 Small Caps % New family
\newfam\scfam
\font\tensc  =cmcsc10

\txtf\scfam=\tensc

\def\sc{\fam\scfam\tensc}

%%%%%%%%%%%%%%%%%%%%%%%%%%%%%%%%%%%%%%%%%%%%%%%%%%%%%%%%%%%%%%%
% 14-, 9-point environments
%%%%%%%%%%%%%%%%%%%%%%%%%%%%%%%%%%%%%%%%%%%%%%%%%%%%%%%%%%%%%%%
% 14-point

\def\frtnmath{%
% \textfont            
\txtf0=\frtnrm        
\txtf1=\frtnmi         
}

\def\frtnpoint{%
\baselineskip=16.8pt plus.5pt minus.5pt%
\def\rm{\fam0 \frtnrm}%
\def\oldstyle{\fam1 \frtnmi}%
\everymath{\frtnmath}%
\everyhbox{\frtnrm}%
\frtnrm }

% 9-point

\def\ninemath{%
% \textfont            
\txtf0=\ninerm        
\txtf1=\ninemi        
\txtf2=\ninesy        
\txtf\itfam=\nineit      
\txtf\bffam=\ninebf      
}

\def\ninepoint{%
\baselineskip=10.8pt plus.1pt minus.1pt%
\def\rm{\fam0 \ninerm}%
\def\oldstyle{\fam1 \ninemi}%
\def\it{\fam\itfam\nineit}%
\def\bf{\fam\bffam\ninebf}%
\everymath{\ninemath}%
\everyhbox{\ninerm}%
\ninerm }

%%%%%%%%%%%%%%%%%%%%%%%%%%%%%%%%%%%%%%%%%%%%%%%%%%%%%%%%%%%%%%%%%%%%%

\def\text#1{\hbox{\rm #1}}

\def\GL{\mathop{\rm GL}\nolimits}

\def\Lie{\mathop{\rm Lie}\nolimits}
\def\Ker{\mathop{\rm ker}\nolimits}
\def\Aut{\mathop{\rm Aut}\nolimits}

\def\IH{\mathop{\hbox{$I$}\!\hbox{$H$}}\nolimits}
\def\IP{\mathop{\hbox{$I$}\!\hbox{$P$}}\nolimits}
\def\Ib{\mathop{\hbox{$I$}\hbox{$b$}}\nolimits}

\def\cite#1{{\uppercase{#1}}}
\def\ref#1{{\uppercase{#1}}}
\def\label#1{{\uppercase{#1}}}

\def\br{\hfill\break} 

%%%%%%%%%%%%%%%%%%%%%%%%%%%%%%%%%%%%%%%%%%%%%%%%%%%%%%%%%%%%%%%%%%%%%
% General Structure
%%%%%%%%%%%%%%%%%%%%%%%%%%%%%%%%%%%%%%%%%%%%%%%%%%%%%%%%%%%%%%%%%%%%%

\def\topmatter{\null\firstpagetrue\vskip\bigskipamount}
\def\endtopmatter{\vskip2\bigskipamount}

\def\title#1{%
\vbox{\raggedright\frtnpoint
\noindent #1\par}
\vskip 2\bigskipamount}        % linebreaks with \br

% \shorttitle =>  short version of the title for the 
% running head

\def\shorttitle#1{\rightheadtext={#1}}              

\newif\ifThanks
\global\Thanksfalse

\def\author#1{\begingroup\raggedright
\noindent{\sc #1\ifThanks$^*$\else\fi}\endgroup
\leftheadtext={#1}\vskip \bigskipamount}

\def\address#1{\begingroup\raggedright
\noindent{#1}\endgroup\vskip\bigskipamount}

\def\endabstract{\endgroup}

\long\def\abstract#1\endabstract{\par
\begingroup\ninepoint\narrower
\noindent{\sc Abstract.\enspace}#1%
\vskip\bigskipamount\endabstract}

\def\section#1#2{\bigbreak\bigskip\begingroup\raggedright
\noindent{\bf #1.\quad #2}\nobreak
\medskip\endgroup\noindent\ignorespaces}

\def\thm#1{\medbreak\noindent{\sc Theorem #1.\enspace}\begingroup
\it\ignorespaces}
\def\endthm{\endgroup\bigbreak}

\def\thmo#1{\medbreak\noindent{\sc #1\enspace}\begingroup
\it\ignorespaces}
\def\endthm{\endgroup\bigbreak}

\def\lem#1{\medbreak\noindent{\sc Lemma #1.\enspace}
\begingroup\it\ignorespaces}
\def\endlem{\endgroup\bigbreak}

\def\state#1{\medbreak\noindent{\sc#1\enspace}
\begingroup\it\ignorespaces}
\def\endstate{\endgroup\bigbreak}  
                  %    for unnumbered statements

\def\qed{$\mathord{\vbox{\hrule\hbox{\vrule
             \hskip5pt\vrule height5pt\vrule}\hrule}}$}

\def\demo{\medbreak\noindent{\it Proof.\enspace}\ignorespaces}
\def\enddemo{\penalty-100\null\hfill\qed\bigbreak}

\newdimen\EZ

\EZ=.5\parindent

\newbox\itembox

\newdimen\ITEM
\newdimen\ITEMORG
\newdimen\ITEMX
\newdimen\BUEXE

\def\iteml#1#2#3{\par\ITEM=#2\EZ\ITEMX=#1\EZ\BUEXE=\ITEM
\advance\BUEXE by-\ITEMX\hangindent\ITEM
\noindent\leavevmode\hskip\ITEM\llap{\hbox to\BUEXE{#3\hfil$\,$}}%
\ignorespaces}

%%%%%%%%%%%%%%%%%%%%%%%%%%%%%%%%%%%%%%%%%%%%%%%%%%%%%%%%%%%%%%%%%%%%%
% Running Heads
%%%%%%%%%%%%%%%%%%%%%%%%%%%%%%%%%%%%%%%%%%%%%%%%%%%%%%%%%%%%%%%%%%%%%

\newif\iffirstpage\newtoks\righthead
\newtoks\lefthead
\newtoks\rightheadtext
\newtoks\leftheadtext
\righthead={\ninepoint\rm\hfill{\the\rightheadtext}\hfill\llap{\folio}}
\lefthead={\ninepoint\rm\rlap{\folio}\hfill{\the\leftheadtext}\hfill}
\headline={\iffirstpage\hfill\else
   \ifodd\pageno\the\righthead\else\the\lefthead\fi\fi}
\footline={\iffirstpage\hfill\global\firstpagefalse\else
     \hfill\fi}

\leftheadtext={}
\rightheadtext={}

%%%%%%%%%%%%%%%%%%%%%%%%%%%%%%%%%%%%%%%%%%%%%%%%%%%%%%%%%%%%%%%%%%%%%
% References
%%%%%%%%%%%%%%%%%%%%%%%%%%%%%%%%%%%%%%%%%%%%%%%%%%%%%%%%%%%%%%%%%%%%%

\def\refs{\bigbreak\bigskip\noindent{\bf References}\medskip
\begingroup\ninepoint}
\def\endrefs{\par\endgroup}

\def\bibiteml#1{\iteml03{[#1]}}%

%%%%%%%%%%%%%%%%%%%%%%%%%%%%%%%%%%%%%%%%%%%%%%%%%%%%%%%%%%%%%%%%%%%%%

\topmatter
\title{On Generalized $h$--Vectors of  \br
Rational Polytopes with a Symmetry of Prime Order}
\shorttitle{On Generalized $h$--Vectors of 
Rational Polytopes with a Symmetry}

\author{Annette A'Campo--Neuen}
\address{%
Fakult\"at f\"ur Mathematik und Informatik\br
Universit\"at Konstanz\br
\br
Annette.ACampo@uni-konstanz.de}
\endtopmatter

\abstract We prove tight lower bounds for the coefficients of the  
generalized $h$-vector of a rational  polytope with a symmetry of 
prime order that is fixed--point--free on the boundary. These bounds 
generalize  results of R.~Stanley and R.~Adin for the $h$--vector 
of a simplicial rational polytope with a central symmetry or a  
symmetry of prime order respectively.
\endabstract

\section{0}{Introduction}
For  simplicial polytopes,  there is a beautiful complete 
characterization of the occuring $h$--vectors, conjectured by 
P.~McMullen and  proved by L.~Billera, C.~Lee, R.~Stanley and  
by McMullen (see [BL], [St1] and [McM]). 
Assuming that in addition the polytope admits a  symmetry, it 
is natural to ask for the resulting restrictions upon the  
corresponding $h$--vector. 

A.~Bj\"orner conjectured tight lower bounds for the coefficients 
of the $h$--vector of a centrally--symmetric simplicial polytope 
that were proved by R.~Stanley in [St2], using the theory of 
toric varieties. After a small perturbation preserving the 
combinatorics as well as the central symmetry one can assume 
that the polytope is rational. Then the polytope defines a 
rationally smooth projective toric variety with an equivariant 
involution and the coefficients of the $h$--vector are the 
Betti--numbers of the  variety. Stanley  in fact proves lower 
bounds for those  Betti--numbers.

Stanley's results were generalized by R.~Adin to the case of 
rational simplicial polytopes admitting a symmetry of prime 
order without fixed points on the boundary of the polytope 
(see [Ad]). Here it is essential to assume  that the polytope 
is rational,  and Adin's statement does not make sense 
without it. 

The aim of this paper is to prove an analogous result to 
Adin's for rational polytopes that are no longer assumed to 
be simplicial.  Then the Betti--numbers of the associated 
projective toric varieties  in general  are not combinatorial 
invariants. We have to replace singular cohomology by 
rational intersection cohomology of middle perversity to 
arrive at combinatorial invariants. The intersection Betti 
numbers can be expressed in terms of numbers of certain 
flags of faces of the polytope and form the so--called 
generalized $h$--vector (see [St3]).  

We obtain lower bounds for the coefficients of the generalized 
$h$--vector of a rational polytope with the same type of 
symmetry as studied by Adin. In his proof Adin uses the 
refined Poincar\'e--series of the Stanley--Reisner--ring 
associated to the polytope. We will generalize this method 
by considering the refined Poincar\'e--series of the 
equivariant intersection cohomology instead. Note that  
in the simplicial case the equivariant intersection 
cohomology forms a ring which is isomorphic to the 
Stanley--Reisner--ring.

In Section~1 we fix the notation, recall Adin's result  and  
state our generalization (see Theorem~1.2). Section~2 is 
devoted to the interpretation of the result in terms of 
projective toric varieties (see Theorem 2.1). Section~3 
contains the facts about equivariant intersection 
cohomology that are needed to complete the proof which 
is carried out in Section~4. 

\section{1}%
%{Convex--Geometrical Background and Statement of the Result}
{Rational Polytopes with Symmetry of Prime Order}
Let $P$ denote a rational polytope in $\RR^n$, i.e.~the 
convex hull of a finite number of points with rational 
coordinates. We assume that $P$ is of full dimension 
and that the origin is placed in its  center of mass. 
By definition, a symmetry of $P$ is a bijection of $P$ 
induced by a  linear map $A$. Note that the assumption 
that $P$ is rational implies  $A\in\GL_n(\QQ)$. 

A symmetry of $P$ does not have fixed points on the 
boundary of the polytope if and only if the inducing 
linear map $A$ does not have $1$ as an eigenvalue. 
This means that the linear map permutes the proper 
faces of the polytope without mapping any face onto 
itself.

We denote the set of proper faces of  $P$ by $\FP$. 
(For technical reasons, we also consider the empty 
set as a face with the convention $\dim\emptyset=-1$.) 
If  $P$ admits a symmetry of prime order $p$ without 
fixed points on the boundary, then all the orbits of 
the induced permutation of  $\FP$  have length $p$, 
except the orbit of the empty face. So in particular,  
the numbers $f_j$ of faces of dimension $j$ of $P$ 
(for $0\leq j\leq n-1$) are all divisible by $p$. 
Moreover,  the dimension $n$ is a multiple of $p-1$, 
since the  linear map $A\in\GL_n(\QQ)$ defining the 
symmetry does not have $1$ as an eigenvalue, and 
therefore its characteristic polynomial  is a power 
of the $p$--th cyclotomic polynomial $1+x+\cdots +x^{p-1}$.

For example, the  $(p-1)$--simplex $S\subset\RR^{p-1}$, 
obtained as the  convex hull of the canonical basis 
vectors $e_1,\dots,e_{p-1}$ and the vector 
$v:=-\sum_{i=1}^{p-1} e_i$, has a symmetry of order $p$, 
induced by  the linear map, sending $e_i$ to $e_{i+1}$ 
for all $i<p-1$ and $e_{p-1}$ to $v$.
Polytopes admitting a symmetry of order $2$ without 
fixed points on the boundary are precisely 
centrally--symmetric polytopes.
 
The  $h$--vector $(h_0,\dots,h_n)$ of an $n$--dimensional 
polytope $P$ is defined as follows:
$$\sum_{j=0}^n h_j x^j:= 
\sum_{j=0}^n f_{j-1} (x-1)^{n-j}\in\ZZ[x] \,,$$
where  $f_j$ denotes the number of $j$--dimensional 
faces of $P$ (for $-1\leq j\leq n-1$). For example, 
$h_n=1$, $h_{n-1}=f_0-n$ and 
$h_0=\sum_{j=0}^n (-1)^{n-j}f_{j-1}$. If $P$ is a simplex, 
then $h_j=1$ for all $j$. Note that the definition of 
the $h$--vector makes sense for arbitrary polytopes. 
But as we mentioned in the introduction, if $P$ is 
rational and simplicial then the coefficients of the 
$h$--vector have a topological interpretation as 
Betti--numbers.

Recall that a polynomial 
$q(x)=\sum_{i=0}^n a_i x^i\in\ZZ[x]$ is called 
{\it symmetric}, if $a_i=a_{n-i}$ for all $i$, and 
it is called {\it unimodal}, if its coefficients 
increase up to a certain index and then decrease again.

In [Ad], R. Adin proved the following theorem, thereby 
generalizing Stanley's result ([St2]) for  
centrally--symmetric simplicial polytopes:

\thm{1.1} {\rm (Adin)}\enspace Let $(h_0,\dots,h_n)$ 
denote the $h$--vector of a rational simplicial 
polytope $P$ with a symmetry of prime order $p$ that 
is fixed--point--free on the boundary. Then $p-1$ 
divides $n$, and the polynomial 
$$\sum_{j=0}^{n} h_j x^j - (1+x+ \cdots +x^{p-1})^{r}
\in\ZZ[x]\,,\quad\text{\it where }r:={n\over p-1},$$ 
is symmetric and unimodal, and
all its coefficients are divisible by $p$.
\endthm

Since $h_0=1$, the constant term of the polynomial 
is zero. So the unimodality of the polynomial  implies 
in particular that the coefficients are nonnegative. 
The resulting lower bounds for the $h_j$ are tight, 
as Adin shows by constructing examples of rational 
simplicial polytopes with a ``minimal'' $h$--vector. 
More precisely, given a natural number $n$ and a 
prime number $p$ such that $p-1$ divides $n$, and 
for $i=1,\dots,r:={n\over p-1}$  a copy $S_i$ of the  
$(p-1)$--simplex $S$ in $V_i:=\RR^{p-1}$, then the 
convex hull of the union $\bigcup_{i=1}^r S_i$ in 
the direct sum  $\bigoplus_{i=1}^r V_i$ has a symmetry 
of order $p$ that is fixed--point--free on the 
boundary. In this case $\sum_{j=0}^{n} h_j x^j=
(1+x+ \cdots +x^{p-1})^{r}$.

In this article we want to further generalize 
Adin's result. Namely, we want to consider   
rational polytopes that are not necessarily simplicial. 
Following Stanley (see [St3]), we  introduce two 
polynomials $h_P$ and  $g_P$ for each polytope $P$, 
that are defined by recursion over the set of faces 
of $P$ as follows:
$$
\leqalignno{
g_{\emptyset} &\equiv 1 \,,&(i)\cr
h_P(x) &=\sum_{F\in\FP} (x-1)^{\dim P-\dim F-1} g_F(x) \,,
&(ii)\cr
g_P(x) &=\tau_{\leq [(\dim Q)/2]}((1-x) h_P(x))\,,
& (iii) \cr
}
$$ 
where $\tau_{\leq r}$ denotes the truncation operator 
$\tau_{\leq r}(\sum_{i=0}^n a_i x^i):=\sum_{i=0}^r a_i x^i$.

The vector formed by the coefficients of the polynomial 
$h_P$  is called the {\it generalized $h$--vector\/} of 
$P$.
Note that $h_P$ depends only on the set of faces  $\FP$ 
as a partially ordered set, and it involves not only the 
face numbers but also the numbers of certain flags of 
faces of $P$. If $P$ is simplicial, then the generalized 
and the ordinary $h$--vector coincide. But this is not 
true for a general rational polytope. For example the 
ordinary $h$--vector of the $3$--dimensional cube  
is $(1,-1,5,1)$ whereas the generalized $h$--vector is 
$(1,5,5,1)$.

The aim of this paper is to prove the following 
generalization of Theorem 1.1:

\thm{1.2} Let $P$ be a rational polytope of dimension $n$ 
admitting a sym\-metry of prime order $p$ without 
fixed points on the boundary, and let 
$(\th_0,\th_1,\dots,\th_n)$ denote its generalized 
$h$--vector. 
Then $p-1$ divides $n$, and the polynomial 
$$\sum_{j=0}^{n} {\tilde h_j} x^j - (1+x+ \cdots +x^{p-1})^{r}\in\ZZ[x]\,,\quad\text{\it where }
r:={n\over p-1}\,, $$
is symmetric and unimodal, and
all its coefficients are divisible by $p$. In particular, 
all the coefficients are nonnegative.
\endthm

\section{2}{Projective Toric Varieties and Intersection 
Cohomology}
Both Stanley and Adin prove their results by first 
translating them to  statements about projective toric 
varieties, and we will follow the same strategy. 
Let $P\subset\RR^n$ be a rational $n$--dimensional 
polytope with a symmetry of prime order $p$, induced 
by a linear map $A\in\GL_n(\QQ)$. Suppose that $A$ 
does not have $1$ as an eigenvalue. Then the 
characteristic polynomial of $A$ equals 
$(1+x+\dots+x^{p-1})^r$ for $r:=n/(p-1)$, 
and as an $A$--module $\QQ^n$ decomposes into a 
direct sum of $r$ copies of  $\QQ[t]/(1+x+\dots+x^{p-1})$. 
So in particular, the rational canonical form of $A$ has 
only integer entries. So after a suitable rational base 
change we can assume that $A\in\GL_n(\ZZ)$. 

The polytope $P$ defines a fan 
$\Delta_P:=\{\RR_{\geq 0}\cdot (F\cup\{0\})\mid F\in\FP\}$ 
in $\RR^n$ consisting of all the cones through proper 
faces of $P$. Since $P$ is rational all the cones are 
rational with respect to the lattice $\ZZ^n$. The pair  
$(\Delta_P,\ZZ^n)$ corresponds to an $n$--dimensional 
projective toric variety $X=X_P$ (see e.g.~[Fu]). 
The vector space $\RR^n$ is naturally identified with 
the Lie--algebra of the maximal compact subgroup $S$ 
of the algebraic torus $T$ acting on $X$, and the 
lattice $\ZZ^n$ in $\RR^n$ is  the kernel $N$ of the 
exponential map from $\Lie(S)$ to $S$.

 Moreover, the linear map $A\in\GL_n(\ZZ)$ defines a 
unique automorphism $\varphi_A$ of  $T$ determined by 
the property  $(d\varphi_A)_e|_N=A$.   
The isomorphism $\varphi_A$ extends to an automorphism 
of the   toric variety $X_P$ since  $A$ permutes the 
cones of $\Delta_P$ (for details on equivariant 
morphisms of toric varieties see e.g.~[Fu]). The 
assumption that $A$ has no fixed points on the 
boundary of $P$ implies that the only 
$\varphi_A$--stable $T$--orbit in $X_P$ is the 
dense orbit.

Conversely, every projective toric variety $X$ with a 
$T$--equivariant automorphism $\varphi$ of prime order 
$p$ can be written in the form $X=X_P$ for some 
polytope $P$ such that the linear map 
$A=d\varphi_e\vert_N\in\GL_n(\ZZ)$ corresponding to 
$\varphi$ induces a symmetry of $P$. If in addition 
$\varphi$ fixes only the dense $T$--orbit then $A$ 
does not have $1$ as an eigenvalue.

If $P$ is simplicial, then the associated projective 
toric variety is rationally smooth, i.e.~it has at 
worst finite quotient singularities. In that case 
the odd Betti numbers of singular cohomology vanish, 
and the even Betti--numbers  are precisely the 
entries of the $h$--vector of $P$, 
i.e.~$b^{2j}(X_P)=h_j$ for $0\leq j\leq n$. So  
the Poincar\'e--polynomial of the toric variety 
$X_P$ is of the form
$$P_X(t):=\sum_{j=0}^{2n} \dim(H^j(X;\QQ)) t^j = 
\sum_{j=0}^{n} h_j t^{2j} \,.$$
By Poincar\'e--duality, the polynomial $P_X$ is 
symmetric, and in particular $h_j=h_{n-j}$ for all 
$j$. These relations for the coefficients of the 
$h$--vector are also known as the Dehn--Sommerville 
equations. 
Moreover, the polynomial $P_X$ is unimodal as a 
consequence of the hard Lefschetz theorem. (The fact 
that the hard Lefschetz theorem for singular 
cohomology is true for rationally smooth varieties 
follows from the fact that there also is a hard 
Lefschetz theorem for intersection cohomology of 
middle perversity of algebraic varieties (see [BBD] 
and [Sa]). If the variety is rationally smooth then 
the intersection cohomology of middle perversity 
and singular cohomology coincide.)  

For example, the $(p-1)$--simplex $S$, for a given 
prime number $p$, corresponds to the toric variety 
$\PP^{p-1}(\CC)$ whose Poincar\'e--polynomial is 
$1+t^2+\cdots+t^{2(p-1)}$, and as stated in Section 1, 
all the entries of the $h$--vector of $S$ are equal 
to $1$. 

If we drop the assumption that $P$ is simplicial then 
the associated projective toric variety $X_P$ can 
have more complicated singularities, and the Betti--numbers 
of  singular  cohomology are no longer combinatorial 
invariants. But instead of singular cohomology, we 
can consider the rational intersection cohomology of 
middle perversity. It turns out that the intersection 
Betti numbers are in fact given by the entries of the 
generalized $h$--vector  $(\th_0,\dots,\th_n)$ of $P$. 
More precisely, the odd intersection Betti numbers 
$\Ib_{2j+1}(X)$ vanish and the even intersection Betti 
numbers satisfy $\Ib_{2j}(X)=\th_j$. (In his survey 
article [St3], Stanley remarks that this result was 
proved independently by J.~Bernstein, A.~Khovanskii 
and R.~MacPherson, but their proofs have not been 
published. 
Proofs appeared in [Fi] by K.~Fieseler and in [DL] 
by J.~Denef and F.~Loeser.) In particular, the 
intersection Betti numbers are combinatorial invariants.

The Poincar\'e--polynomial of the intersection 
cohomology of middle perversity 
$$\IP_X(t):=\sum_{j=0}^{2n} \dim(\IH^j(X;\QQ)) 
t^j=\sum_{j=0}^{2n} \Ib_{j} t^j$$
 is symmetric, since Poincar\'e--duality holds for 
intersection cohomology.
Moreover, as mentioned above there is also a hard 
Lefschetz theorem for intersection cohomology of 
middle perversity  and therefore $\IP_X$ is unimodal 
(see [St3], Cor.~3.2).

So in terms of  toric varieties Theorem 1.2  reads 
as follows:
\thm{2.1} Let $X$ be a projective toric variety  of 
dimension $n$, and suppose that $X$ has a 
$T$--equivariant automorphism $\varphi$ of  prime 
order $p$, defining a fixed--point--free permutation 
of the $T$--orbits in the boundary of $T$ in $X$. 
Then $p-1$ divides $n$, and the polynomial 
$$q(x):=\sum_{j=0}^{n} \Ib_{2j} x^j - (1+x+ \cdots +x^{p-1})^{r}\in\ZZ[x]\,,\quad\text{\it where }
r:={n\over p-1}\,, $$
is symmetric and unimodal, and
all its coefficients are divisible by $p$. 
\endthm

The symmetry of the polynomial $q$ immediately follows 
from the symmetry of the Poincar\'e--polynomial $\IP_X$.
 So to prove the theorem, we only have to show that $q$ 
is unimodal and that the coefficients are divisible by 
$p$.
Analogous to Adin's proof of Theorem 1.1, we will 
proceed by interpreting ${1\over p}q$ as a  polynomial 
whose coefficients count certain dimensions. 
In our case, these are the dimensions of the 
$\chi$--eigenspaces of $\varphi$ on the graded pieces 
of the rational intersection cohomology of $X$ for 
some primitive $p$--th root of unity $\chi$.

The automorphism $\varphi$ induces a linear map on the 
$j$--th 
intersection cohomology $\IH^j(X;\QQ)$ of $X$ for 
$0\leq j\leq 2n$. Since $\varphi$ is of order $p$, the 
possible eigenvalues of the induced endomorphism are 
the elements of the group $G$ of $p$--th roots of unity. 
So over the complex numbers, we have the following 
decomposition:
$$\IH^j(X;\QQ)\otimes\CC=\bigoplus_{\chi\in G} 
\IH^j_\chi(X)\,, $$
where $\IH^j_\chi(X)$ denotes the $\chi$--eigenspace of 
$\varphi$.
The crucial result is the following proposition.

\thmo{Proposition 2.2} Let $X$ be a projective toric 
variety  of dimension $n$, and suppose that $X$ has a 
$T$--equivariant automorphism $\varphi$ of  prime order 
$p$, defining a fixed--point--free permutation of the 
$T$--orbits in the boundary of $T$ in $X$. Let 
$r:={n\over p-1}$. Then 
$$P_{\chi}(t):=\sum_{j=0}^{2n} (\dim\IH^j_\chi(X))  
t^j = {1\over p}  
(\IP_X(t) - (1+t^2+\cdots + t^{2(p-1)})^r) $$
for every primitive $p$--th root of unity $\chi$.

In particular, $\IP_X(t) = (1+t^2+\cdots + t^{2(p-1)})^r$ 
iff $\varphi$ induces the identity on $\IH^*(X;\QQ)$.
\endthm

The proof of this proposition will be given in the last 
section using equivariant intersection cohomology. We 
conclude this section by deducing Theorem 2.1 from the 
proposition.
\bigskip

\noindent{\it Proof of Theorem 2.1.\enspace}  It suffices 
to show that $P_{\chi}$ is unimodal. The intersection 
cohomology of middle perversity $\IH^*(X;\QQ)$  forms 
a module over the ring of singular cohomology $H^*(X;\QQ)$. 
The hard Lefschetz theorem for intersection cohomology 
asserts that there is an element $\omega\in H^2(X;\QQ)$ 
corresponding to a generic hyperplane section of some 
embedding of $X$ into complex projective space, such 
that the multiplication with $\omega$ induces an 
injective map from $\IH^{2j}(X;\QQ)$ to 
$\IH^{2j+2}(X;\QQ)$ for all $0\leq j\leq [{n\over 2}]$. 

The cohomology class $\omega$ is $\varphi$--invariant 
since the action of $\varphi$ on $X$ extends to a 
linear action on the projective space (see the argument 
given in [St2], provided by S.~Kleiman). Therefore 
the multiplication with $\omega$ commutes with the 
action of $\varphi$ on the intersection cohomology, 
and hence $\omega\cdot\IH^{2j}_\chi(X)\subset
\IH^{2j+2}_\chi(X)$ for all $0\leq j\leq [{n\over 2}]$ 
and for all characters $\chi$. So in particular, all 
the polynomials $P_{\chi}$ are unimodal.
\enddemo

\section{3}{Equivariant Intersection Cohomology}
In this section we first briefly recall some definitions, 
and then we state the facts that are used in the proof 
of Proposition 2.2.
For the algebraic torus $T=(\CC^*)^n$, the principal 
$T$--bundle $ET:=(\CC^{\infty}\backslash\{0\})^n\to 
(\PP^{\infty}(\CC))^n=:BT$ is a classifying bundle 
since the action of $T$  on the contractible space 
$ET$ by componentwise scalar multiplication is free. 
The cohomology ring $H^*(BT;\QQ)$ of the base space 
$BT$ is a polynomial ring over $\QQ$ in $n$ variables, 
where each variable has degree $2$, and the 
corresponding Poincar\'e--series is 
$(1+t^2+t^4+\cdots)^n=(1-t^2)^{-n}$.

Now let $X$ be a toric variety acted on by the torus 
$T$, and set $E_q:=(\CC^{q+1}\backslash\{0\})^n$ for 
$q\in\NN\cup{\infty}$. The diagonal action of $T$ on 
$E_q\times X$ is free, and forming the quotient by 
$T$ we obtain  $X_{T,q}:= (E_q\times X)/T= E_q\times_T X$.  
The natural projection map 
$p_q:X_{T,q}\to (\PP^q(\CC))^n=:B_q$ is a fibration 
with fiber $X$. 

The intersection cohomology group $\IH^j(X_{T,q};\QQ)$  
is independent of $q$ for $0\leq j<2q$. Therefore it 
is natural to define equivariant intersection 
cohomology of $X$ with rational coefficients in the 
following way (see (2.12) in [Ki]):
$$\IH^j_T(X;\QQ):= \IH^j(X_{T,q};\QQ) \quad\text{for 
$j\geq 0$, \ where $q>{j/2}$.}$$
   
The equivariant intersection cohomology of $X$ is 
related to its ordinary intersection cohomology  
in the following way (see e.g. (2.14) in [Ki]):

\state{Remark.}
Let $X$ be a projective toric variety. Then
$$\IH^*_T(X;\QQ) \simeq H^*(BT;\QQ) \otimes_{\QQ} 
\IH^*(X;\QQ).
\leqno{(3.1)}$$
\endstate

This remark can be proved by observing that the 
spectral sequences associated to the fibrations 
$p_q$ degenerate at the $E_2$--level. So for every $j$, 
there is a filtration of $\IH^j(X_{T,q};\QQ)$ such 
that the occuring factors are precisely the 
$E_2$--terms of the spectral sequence on the 
$j$--th diagonal. On the other hand the direct sum 
of all $E_2$--terms equals the tensor product 
$H^*(B_q;\QQ) \otimes_{\QQ} \IH^*(X;\QQ)$.  

For the  Poincar\'e--series of the equivariant 
intersection cohomology
$\IP_X^T(t):=\sum_{j=0}^{\infty} 
(\dim\IH^j_T(X;\QQ)) t^j$
we conclude from (3.1):
$$
\IP_X^T(t)= (1-t^2)^{-n} \IP_X(t)\,.
\leqno{(3.2)}$$

We will also use the following result, proved 
by K.~Fieseler in [Fi]: 

\thmo{Proposition 3.3.}%\label{Fieseler} 
{\rm (Fieseler)} Let $X$ be a projective toric 
variety corresponding to the fan $\Delta\in\RR^n$. 
Then 
$$\IH^*_T(X;\QQ)\simeq\bigoplus_{\sigma\in\Delta}
\IH^*_T(U_{\sigma},U_{\sigma}\backslash B_{\sigma};\QQ)\,$$
where for every cone $\sigma\in\Delta$, 
$U_{\sigma}$ denotes the corresponding affine chart 
of $X$  and $B_{\sigma}$ denotes the corresponding 
$T$--orbit.
\endthm

\section{4}{Proof of Proposition 2.2}
Returning to the setting of Theorem 2.1 and 
Proposition 2.2, let  us assume that $X$ is a 
projective toric variety of dimension $n$ with 
a $T$--equivariant automorphism $\varphi$ of  
order $p$, for some prime number $p$, and suppose 
that the only $\varphi$--stable $T$--orbit in 
$X$ is the dense orbit.  
Then as remarked in Section 2, the automorphism 
$A:=d\varphi_e\in\Aut(N)$ of the lattice 
$N=\Ker\exp$ does not have $1$ as an eigenvalue. 
Since $A$ is of order $p$, the characteristic 
polynomial of $A$ is a power of the $p$--th 
cyclotomic polynomial.

We will prove Proposition 2.2 by analysing the 
induced action of $\varphi$ on the terms 
occuring in (3.1), and considering their 
respective refined Poincar\'e--series. 

The {\it refined\/} Poincar\'e--series 
$P_W^{\varphi}$ of a graded complex vector 
space $W=\bigoplus_{j=0}^{\infty} W^j$ with 
respect to the action of some graded 
endomorphism $\varphi$ of order $p$ is defined 
as the following formal power series in one 
variable over the group ring $\ZZ[G]$ of the 
group $G$ of $p$--th roots of unity:
$$P_W^{\varphi}(t):=\sum_{j=0}^{\infty} 
\sum_{\chi\in G} (\dim W_{\chi}^j) \,\chi\,t^j  \,,$$
where $W_{\chi}^j$ denotes the eigenspace of 
$\varphi$ in $W^j$ for the eigenvalue $\chi\in G$. 
If for example $W$ is a polynomial ring in one 
variable of degree $d$ which is an eigenvector of 
$\varphi$ for the eigenvalue $\lambda$,  then 
$P_{W}^\varphi(t)=
1+\lambda t^d + \lambda^2 t^{2d}+ \cdots=
(1-\lambda t^d)^{-1}$.

For $W=\IH^*(X;\QQ)\otimes \CC$ or 
$\IH_T^*(X;\QQ)\otimes \CC$,   we denote 
the corresponding refined Poincar\'e--series 
by $\IP_{X}^\varphi$ or $\IP_{X}^{T,\varphi}$ 
respectively. We claim that:
$$\IP_X^{T,\varphi}(t)=P_{BT}^\varphi(t) 
\IP_X^{\varphi}(t) \,.\leqno{(4.1)}$$ 

To justify (4.1) note that the $T$--equivariant 
automorphism $\varphi$ of $X$  also determines  
automorphisms  of $E_q$, as well as of $X_{T,q}$ 
and $B_{q}$, that commute with the fibrations  
$p_q: X_{T,q} \to B_q$. So $\varphi$  acts 
naturally on all the terms occuring in the 
associated spectral sequences. The proof of 
(3.1) implies that  for every $j$ there is an 
isomorphism between the associated graded 
$\varphi$--modules of $\IH^j_T(X;\QQ)$ and 
of $(H^*(BT;\QQ) \otimes_{\QQ} \IH^*(X;\QQ))^j$. 
Since $\varphi$ is of finite order and hence 
semi--simple, we can even conclude that 
$\IH^j_T(X;\QQ)$ and $(H^*(BT;\QQ) \otimes_{\QQ} 
\IH^*(X;\QQ))^j$ are isomorphic as 
$\varphi$--modules and the claim follows.

Let us now compute $P_{BT}^\varphi$. There is a 
natural identification of $H^2(BT;\RR)$ with the 
dual of the Lie--algebra $\Lie(S)$ of the maximal 
compact torus $S\subset T$. On the other hand, 
$\Lie(S)$ is generated as a real vector space  
by $N=\ker\exp$ in $\Lie(S)$. So altogether we 
have a natural isomorphism between $H^2(BT;\QQ)$ 
and $N_{\QQ}^\vee$, and we can view $H^*(BT;\QQ)$ 
as the symmetric algebra of $N_{\QQ}^\vee$.

The action of $\varphi$ on 
$H^2(BT;\QQ)\simeq N_{\QQ}^\vee$ is given by  
the dual of $A$. So the eigenvalues are precisely 
the primitive $p$--th roots of unity, and each 
$\chi\in G\backslash\{1\}$ occurs with the same 
multiplicity $r:=n/(p-1)$. Over the complex 
numbers we can choose a basis of the second 
cohomology consisting of eigenvectors for 
$\varphi$, and  the cohomology ring of $BT$ 
is a polynomial ring in those basis vectors. 
Each $\chi$--eigenvector in the basis of the 
second cohomology contributes a factor 
$(1-\chi t^2)^{-1}$ to the refined 
Poincar\'e--series of $BT$, and hence
$$
P_{BT}^\varphi(t)=\sum_{j=0}^{\infty} \sum_{\chi\in G} 
(\dim H^j_{\chi}(BT)) \chi t^j  = \prod_{1\not=
\chi\in G} {1\over (1-\chi t^2)^{r}}  \,.
\leqno{(4.2)}$$

To describe the  action of $\varphi$ on 
$\IH^*_T(X;\QQ)$ we use that 
by Proposition 3.3, 
$\IH^*_T(X;\QQ)\simeq\bigoplus_{\sigma\in\Delta} \IH^*_T(U_{\sigma},U_{\sigma}\backslash B_{\sigma};\QQ)$, 
where $\Delta$ denotes the fan associated to the 
toric variety $X$. That means that the long exact 
$\IH^*_T$--sequences associated to the pairs 
$(V_r,V_{r-1})$, where 
$V_r:=\bigcup_{\dim\sigma\leq r} B_{\sigma}$, 
split into short exact sequences. Since the 
$\varphi$--action preserves the open exhaustion 
of $X$ by the open subsets $V_r$, $\varphi$ 
acts naturally on all the terms occuring in 
those sequences. Using that $\varphi$ is of 
finite order, we obtain as above that 
$\IH^j_T(X;\QQ)$ and $\bigoplus_{\sigma\in\Delta} \IH^j_T(U_{\sigma},U_{\sigma}\backslash B_{\sigma};\QQ)$ 
are isomorphic even as $\varphi$--modules. So in 
particular, the refined Poincar\'e--series of both 
sides of (3.3) are equal.

The action of $\varphi$ on the right hand side of (3.3) 
permutes the direct summands in the same way as the 
linear map $A$ considered as an element of $\GL_n(\RR)$ 
permutes the cones of $\Delta$.
The only $A$--stable cone of $\Delta$  is the zero--cone. 
We have $B_{0}=T=U_{0}$ and therefore 
$\IH^*_T(U_{0},U_{0}\backslash B_{0};\QQ)=\IH^*_T(T;\QQ)$, 
which consists of one copy of $\QQ$ in degree~$0$. 
So the zero--cone only contributes a one--dimensional 
subspace in the eigenspace for the eigenvalue $1$  
of $\IH^0_T(X;\QQ)$.

Every cone $\sigma\not=0$ has an orbit of length $p$ 
under $A$, and for every $j$ the $\varphi$--stable 
subspace 
$\bigoplus_{k=1}^p 
\IH^j_T(U_{A^k(\sigma)},U_{A^k(\sigma)}\backslash 
B_{A^k(\sigma)};\QQ)$ of $\IH_T^j(X;\QQ)$ decomposes 
into a direct sum of copies of the regular 
representation of the group generated by $\varphi$. 
That implies that in $\bigoplus_{0\not=\sigma\in\Delta} \IH^j_T(U_{\sigma},U_{\sigma}\backslash B_{\sigma};\QQ)$ 
every eigenvalue $\chi\in G$ appears with the same 
multiplicity. 

We obtain:
$$\IP^{T,\varphi}_{X}(t)-1 = {1\over p}\sum_{\chi\in G}  
(\IP_{X}^T(t)-1) \chi  \,.$$
Setting $\mu:={1\over p} \sum_{\chi\in G} \chi \in \ZZ[G]$ 
and using (3.2), we arrive at
$$ 
\IP_X^{T,\varphi}(t)={\mu\over (1-t^2)^n} \IP_X(t) + (1-\mu) \,.
\leqno{(4.3)}$$

Inserting (4.2) and (4.3) into (4.1) we get:
$$
\IP_X^{\varphi}(t)=
\left({\mu\over (1-t^2)^n} \IP_X(t) + (1-\mu)\right) 
\prod_{1\not=\chi\in G} (1-\chi t^2)^{r} \,.
\leqno{(4.4)}$$

\lem{4.5}{\rm In the polynomial ring $\ZZ[G][x]$ over 
the group ring $\ZZ[G]$ of  the group $G$ of $p$--th 
roots of unity the following identities hold:  
\item{(i)} $\mu \prod_{1\not=\chi\in G} (1-\chi x)^{r}=
\mu (1-x)^n$; 
\item{(ii)} $(1-\mu)\prod_{1\not=\chi\in G} (1-\chi x)^r=
(1-\mu)(1+x + \cdots + x^{p-1})^r$,\par
where as above  $\mu={1\over p} \sum_{\chi\in G}\chi$.
}
\endlem

\demo (i) It suffices to check $\mu\chi=\mu$ for all 
$\chi\in G$. But this is clear since the multiplication 
of $\sum_{\chi'\in G} \chi'$ with an element $\chi$ of 
$G$ only permutes the summands.

(ii)  This part of the lemma is contained in [Ad] as 
Lemma 13.
\enddemo

Applying  Lemma 4.5 and inserting $t^2$ for $x$, we can 
rewrite equation (4.4) in the following form:
$$\IP_X^{\varphi}(t)=\mu \IP_X(t) + (1-\mu) 
(1+t^2 + \cdots + t^{2(p-1)})^r \,,$$ 
which is equivalent to:
$$
\IP_X^{\varphi}(t)= (1+t^2+\cdots + t^{2(p-1)})^r + 
{1\over p} \sum_{\chi\in G} (\IP_X(t) - 
(1+t^2+\cdots + t^{2(p-1)})^r)\chi \,.
$$

Now  Proposition 2.2 is an immediate consequence.

\bigskip
\noindent{\bf Acknowledgements.\enspace} I would like 
to thank K.~Fieseler for valuable discussions and the 
referees for useful remarks and suggested improvements.

\refs

\bibiteml{Ad} R.~M.~Adin, {\it On face numbers of 
rational simplicial polytopes with symmetry}, 
Adv. Math.~{\bf 115} (1995), 269--285.

\bibiteml{BBD} A.~Beilinson, J.~Bernstein, P.~Deligne, 
{\it Faisceaux pervers}, 
Ast\'erisque {\bf 100} (1982), 5--171.

\bibiteml{BL} L.~J.~Billera and C.~W.~Lee, 
{\it A proof of the sufficiency of McMullen's 
condition for $f$--vectors of simplicial 
convex polytopes}, J. Combin. Theory (A) {\bf 31} 
(1981), 237--255.

\bibiteml{DL} J.~Denef and F.~Loeser, 
{\it Weights of exponential sums, intersection 
cohomology, and Newton polyhedra}, 
Inv.~Math.~{\bf 106} (1991), 275--294.

\bibiteml{Fi} K.--H.~Fieseler, {\it Rational 
intersection cohomology of projective toric 
varieties}, J. reine angew.~Math.~{\bf 413} 
(1991), 88--98.

\bibiteml{Fu} W.~Fulton, Introduction to Toric Varieties, 
Annals of Math. Studies {\bf 131}, Princeton University 
Press 1993.

\bibiteml{Ki} F.~Kirwan, {\it Intersection homology and 
torus actions}, 
J.~Amer. Math.~Soc.~{\bf 1} (1988), 385--400.

\bibiteml{McM} P.~McMullen, {\it The polytope algebra}, 
Adv. in Math. {\bf 87} (1989), 76--130.

\bibiteml{Sa} M.~Saito, {\it Mixed Hodge modules}, 
Publ. Res. Inst. Math. Sci. (Kyoto Univ.) {\bf 26} 
(1990), 221--333.

\bibiteml{St1} R.~Stanley, {\it The number of faces 
of a simplicial convex polytope}, Adv. Math.~{\bf 35} 
(1980), 236--238.

\bibiteml{St2} R.~Stanley, {\it On the number of faces 
of centrally--symmetric simplicial polytopes}, 
Graphs Combin.~{\bf 3} (1987), 55--66.

\bibiteml{St3}  R.~Stanley, {\it Generalized $h$--vectors, 
intersection cohomology of toric varieties, and related 
results}, in: Commutative Algebra and Combinatorics 
(M.~Nagata and H. Matsumura, eds.), Advanced Studies 
in Pure Math..~{\bf 11}, Kinokuniya, Tokyo, and 
North--Holland, Amsterdam/New York, 1987, pp.~187--213. 

\bibiteml{St4} R.~Stanley, Combinatorics and 
Commutative Algebra, Progress in Math.~{\bf 41}, 
Birk\-h\"auser Boston etc., Second Edition 1996.

\endrefs
\bye